\documentclass[a4paper,10pt,fleqn]{article}
\usepackage{graphics,graphicx}
\usepackage{natbib}
\usepackage{amsmath}
\usepackage{amssymb}

\chardef\at=`\@
\DeclareRobustCommand{\qed}{%
  \ifmmode 
  \else \leavevmode\unskip\penalty9999 \hbox{}\nobreak\hfill
  \fi
  \quad\hbox{\qedsymbol}}
\newcommand{\mathbold}[1]{\mbox{\boldmath $#1$}}

\newcommand{\pr}{\mathbold P}

\newcommand{\argmin}{\operatorname{argmin}}

\begin{document}
\begin{center}
{\large On stepwise regression\\}
\quad\\
Laurie Davies\\
Fakul\"at Mathematik\\
Universit\"at Duisburg-Essen\\
\texttt{laurie.davies@uni-due.de}
\end{center}
\quad\\
\begin{abstract}
Given data $y$ and $k$ covariates $x$ one problem in linear regression
is to decide which in any of the covariates to include when regressing
$y$ on the $x$. If $k$ is small it is possible to evaluate each subset
of the $x$. If however $k$ is large then some other procedure must be
use. Stepwise regression and the lasso are two such procedures but
they both assume a linear model with error term. A different approach
is taken here which does not assume a model. A covariate is included if
it is better than random noise. This defines a procedure which is
simple both conceptually and algorithmically.
\end{abstract}

\section{Introduction} \label{sec:intro}
In a forward stepwise regression the next variable to be included is
the one which gives the largest reduction in the sum of squared
residuals. The decision as to whether to include this variable is
based on the result of an $F$-test which in turn assumes a linear
model
\begin{equation}  \label{equ:st_model}
Y=x^t\beta+\varepsilon.
\end{equation}
The $f$-test does not take the adaptive nature of the procedure into
account. Such a test exist for the lasso  (\cite{LOKTAYTIB214}) but it
is also based on the model (\ref{equ:st_model}) and requires
assumptions for its validity.

The procedure described below is based on \cite{DAV16} . It does not assume a model
 and consequently makes no assumptions about about an error term or
 about the covariates. It is in other words a procedure. It is
 legitimate and  possible to investigate its behaviour under the model
 (\ref{equ:st_model}) but this will not be done here.

\section{The procedure}

\subsection{Least squares}
Suppose that of the $k$ covariates $k_1<k$ of them have already been
included and their sum of squares by $ss(k_1)$. There remain
$k_0=k-k_1$ covariates and the candidate for 
inclusion is the one whose inclusion decreases the sum of squared
residuals by the most.  Denote this sum of squared residuals by
$ss(k_0)$ so that the reduction in the sum of squares is
\begin{equation}
ss(k_1)-ss(k_0).
\end{equation}

Replace nor the $k_0$ covariates by i.i.d. $N(0.1)$ random
variables. If one of these is included together with the $k_1$
covariates already included it is a simple exercise to see that the
sum of squared residuals is approximately
\begin{equation} \label{equ:chi2}
ss(k_1)-\frac{ss(k_1)}{n}\chi^2_1 .
\end{equation}
Choosing  those random variable which lead to the largest reduction in
the sum of squares leads to a reduction
\begin{equation} \label{equ:chi2_m}
ss(k_1)-\frac{ss(k_1)}{n}\max\{\chi^2_1,
\ldots,\chi^2_1\}=\frac{ss(k_1)}{n}E(k_0)
\end{equation}
where the maximum is taken over $k_0$ independent $\chi^2_1$ random
variables. The probability that the best of  random variables is
better than the best of the remaining remaining $k_0$
covariates is therefore
\begin{equation}
\pr\left(ss(k1)-\frac{ss(k_1)}{n}E(k_0) <ss(k_1)-ss(k_0)\right)
\end{equation}
or equivalently
\begin{equation} 
\pr\left(E(k_0)> \frac{n}{ss(k1)}\left(1-\frac{ss(k_0)}{ss(k_1)}\right)\right).
\end{equation}
If this probability is reasonably large, say 0.1, then in $10\%$ of the
cases the included covariate is no better than random noise. This
probability must be specified in advance by a number $\alpha$. If 
\begin{equation}
\pr\left(E(k_0)>
  \frac{n}{ss(k1)}\left(1-\frac{ss(k_0)}{ss(k_1)}\right)\right)< \alpha
\end{equation}
then the covariate is included. Otherwise the procedure is
terminated.

As $E(k_0)$ is the maximum of $k_0$ $\chi^2_1$ random variables  
\begin{equation} \label{equ:dist_E}
\pr(E(k_0)>x)=1-\text{pchisq}(x,1)^{k_0}
\end{equation}
where $\text{pchisq}(x,\nu)$ is the distribution function of a  $\chi^2$
random variable with $\nu$ degrees of freedom. The covariate is
therefore included if
\begin{equation}
 \frac{n}{ss(k1)}\left(1-\frac{ss(k_0)}{ss(k_1)}\right)>\text{qchisq}((1-\alpha)^{1/k_0},1).
\end{equation}
where $\text{qchisq}(x,\nu)$ is the inverse distribution function of a  $\chi^2$
random variable with $\nu$ degrees of freedom. More informatively one
can calculate the $P$-value
\begin{equation} \label{equ:pval}
1-\text{pchisq}\left( \frac{n}{ss(k1)}\left(1-\frac{ss(k_0)}{ss(k_1)}\right)\right)^{k_0}
\end{equation}

\subsection{$M$-regression}
The method can in principle (with the obvious modifications) be
applied to $L_1$ regression but with the disadvantage that there does
not exist a simple expression corresponding to (\ref{equ:chi2}). If
there is a particular interest in $L_1$ regression then simulations
will be required. If however $L_1$ regression is only used as a
protection against outlying $y$-values this can also be provided by
$M$-regression for which a version of (\ref{equ:chi2}) is available.

Let $\rho$ by a symmetric positive twice differentiable convex function
with $\rho(0)=0$.  The default function will be the $\rho$ function
used in \cite{DAV14}, namely
\begin{equation}
\rho_c(u)=\left\{\begin{array}{ll}
\vert u\vert, & \vert cu\vert \ge 15\\
2\log(0.5+0.5\exp(cu))/c-u,&\vert cu\vert <15\\
\end{array}
\right.
\end{equation}
where $c$ is a tuning constant with default value $c=1$. An
alternative choice could be Huber's $\rho$-function with a tuning
constant (\cite{HUBRON09}). The sum of squared residuals $ss(k)$ is replaced
by
\begin{equation} 
s_{\rho}(k)=\argmin_{\beta}\sum_{i=1}^k \rho(y_i-x_i^t\beta) .
\end{equation}

As it stands $s_{\rho}(k)$ is not satisfactory and must be augmented
by a data dependent scale value $\sigma$ to give
\begin{equation} 
s_{\rho}(k,\sigma)=\argmin_{\beta}\sum_{i=1}^k \rho\left(\frac{y_i-x_i^t\beta}{\sigma}\right) .
\end{equation}
For a given $\rho$ and $\sigma$ $s_{\rho}(k,\sigma)$ can be calculated
using the algorithm described in {\bf 7.8.2} of
\cite{HUBRON09}. Typically only a few number of iterations are
required. 

Given all this (\ref{equ:chi2}) is replaced by
\begin{equation}
s_{\rho}(k_1,\sigma)-\frac{\sum_{i=1}^n
  \rho^{(1)}\left(\frac{r_i}{\sigma}\right)^2}{\sum_{i=1}^n
  \rho^{(2)}\left(\frac{r_i}{\sigma}\right)}\chi^2_1=s_{\rho}(k_1,\sigma)-\frac{s_{\rho^{(1)}}(k_1,\sigma)}{s_{\rho^{(2)}}(k_1,\sigma)}\chi^2_1
\end{equation}
where $\rho^{(1)}$ and $\rho^{(2)}$ are th first and second
derivatives of $\rho$ respectively,
\[r_i=y_i-x_i^t{\hat \beta} \quad \text{with}\quad {\hat
  \beta}=\argmin_{\beta}\sum_{i=1}^{k_1}
\rho\left(\frac{y_i-x_i^t\beta}{\sigma}\right)\] 
and
\[s_{\rho^{(1)}}(k_1,\sigma)=\sum_{i=1}^n
  \rho^{(1)}\left(\frac{r_i}{\sigma}\right)^2, \quad s_{\rho^{(2)}}(k_1,\sigma)=\sum_{i=1}^n
  \rho^{(2)}\left(\frac{r_i}{\sigma}\right) .\]
The $P$-value  (\ref{equ:pval}) becomes
\begin{equation} \label{equ:pval_m}
1-\text{pchisq}\left( \frac{s_{\rho^{(1)}}(k_1,\sigma)}{s_{\rho^{(2)}}(k_1,\sigma)}\left(1-\frac{s_{\rho}(k_0,\sigma)}{s_{\rho}(k_1,\sigma)}\right)\right)^{k_0}
\end{equation}

It remains to specify the choice of scale $\sigma$. The procedure
described here uses a $\sigma$ dependent on the $k_1$ variables
already incorporated. This same $\sigma$ is used to judge whether a
new variable is to be included. This is why there is only one value of
$\sigma$ in (\ref{equ:pval_m}). One possibility is to do a full
$M$-regression and for both location and scale based on the $k_1$
covariates and take $\sigma$ to be the scale part
(\cite{HUBRON09}). This has a certain intellectual consistency but at
the expense of greater programming effort. Instead the following
procedure will be used. If a new covariate is to be included then the
residuals $r_i$ are calculated from an $M$-regression using the $k_1+1$
covariates but based on the $\sigma$ for the original $k_1$
covariates. The new $\sigma$ is taken to be the median absolute
deviation of the  $r_i$ multiplied by the Fisher consistency  factor
1.48 which is the default version of the MAD in {\texttt R}. The
procedure is started using the residuals from best $L_1$ single
covariate calculated using for example \cite{KOEN10}.

\section{Two examples}
The method will be illustrated using the prostate cancer data also
used in (\cite{LOKTAYTIB214}) and the low birth weight data from
\cite{HOSLEM89}. 

The prostate cancer data were obtained from
\cite{LOKVENTUR14}. They are described in \cite{HASTIBFRI08}. The
sample size is $n=97$ with eight covariates. Table~\ref{tab:prostate}
gives the order in which the covariates entered the regression
together with their $P$-values for the $L_2$ and $M$ regressions. The
order was the same for both.   Table~\ref{tab:prostate_1} is the same
but with the first $y$ value changed from  -0.4307829 to 10. It shows
that the results for the $M$ regression remain stable but those for
the $L_2$ regression change considerably apart from the covariate {\it
  lcavol}.

\begin{table}[!]
\begin{center}
\begin{tabular}{|l|c|c|}
\hline
&\multicolumn{2}{c|}{$P$-value}\\
\cline{2-3}
covariate&$L_2$&$M$\\
\hline
lcavol& 0.0000&0.0000\\
lweight&0.0122&0.0083\\
 svi&0.0123&0.0101\\
 lbph&0.4233&0.3408\\
age&0.4952&0.4083\\
pgg45&0.5541&0.4839\\
lcp&0.4093&0.2845\\
gleason&0.7636& 0.7300\\
\hline
\end{tabular}
\caption{The prostate data: the covariates in order of inclusion and their
  $P$-values\label{tab:prostate}}
\end{center}
\end{table}

\begin{table}[b]
\begin{center}
\begin{tabular}{|l|c||l|c|}
\hline
\multicolumn{2}{|c||}{$L_2$}&\multicolumn{2}{|c|}{$M$}\\
\hline
covariate&$P$-value&covariate&$P$-value\\
\hline
lcavol& 0.0000&lcavol&0.0000\\
svi&0.1234&svi&0.0176\\
age&0.6623&lweight&0.0366\\
 lbph&0.4534&lbph&0.4676\\
lweight&0.950&age&0.1766\\
pgg45&0.7615&pgg45&0.5309\\
lcp&0.7615&lcp&0.3337\\
gleason&0.8949&gleason&0.8269\\
\hline
\end{tabular}
\caption{The prostate data but with  $y(1)=10$: the covariates in order of inclusion and their
  $P$-values\label{tab:prostate_1}}
\end{center}
\end{table}

The dependent variable in the low birth weight data is taken to be the
weight of the child. The covariates are:\\
(1) Age of mother, (2) Weight of mother, (3) Smoking status, (4)
History of premature labor, (5) History of hypertension, (6) Uterine
irritability, (7) Number of physician visits, (8) Race-1, (9)
Race-2.\\
Model and functional choice for this data set has been considered in
\cite{DAV14} and \cite{HJOCLA03B} (model choice) and \cite{DAV16}
(functional choice). Table~\ref{tab:lowbirth} gives the
results for the stepwise functional choice. The oder of the covariates
is the same for both methods. The choice (6,9,3) with $\alpha=0.05$
corresponds to the functional encoded as 292 in \cite{DAV16} which one
of the functionals chosen after considering all subsets with
$\alpha=0.05$. The subset (6,9,3,5) corresponds to the functional
encoded as 308 with $\alpha=0.1$ in  \cite{DAV16}.

\begin{table}[b]
\begin{center}
\begin{tabular}{|c|c|c|}
\hline
&\multicolumn{2}{c|}{$P$-value}\\
\cline{2-3}
covariate&$L_2$&$M$\\
\hline
6& 0.0009&0.0008\\
9&0.0187&0.0223\\
3& 0.0015&0.0009\\
5&0.0934&0.1017\\
2& 0.0778&0.0649\\
8&0.8842&0.8616\\
4& 0.9285&0.9038\\
1& 0.8779&0.8359\\
7&0.7557&0.7607\\
\hline
\end{tabular}
\caption{The low birth weight data: the covariates in order of inclusion and their
  $P$-values.\label{tab:lowbirth}}
\end{center}
\end{table}

\end{document}